# CONVERGENCE OF SOME PERTURBED SEQUENCES OF RATIONAL POWERS AND APPLICATION TO SYRACUSE PROBLEM

Hassan Douzi*

*Ibn Zohr University, Faculty of Science, Department of Mathematics, BP8106 Agadir, Morocco.
Tel: 212 6 61224366, Mail: h.douzi @ uiz.ac.ma

## Abstract

*Sequences of rational powers $(\xi(p/q)^n)_{n \geq 0}$, especially in the case $p/q = 3/2$, have a connection with many important combinatorics and number theory problems as for example Syracuse, Z-number and waring problems. Conjectures from such problems are known to be intractable and only few partial results exist until now. In this paper, we study a family of perturbed sequences of rational powers called "Branch sequences" under the form $(S_n = (\xi + \Sigma_n)(p^n/q^{n+e_n}))_{n \geq 0}$. Under the assumption that such sequences are deterministic and they have controlled positive perturbations, we establish the convergence result: $\min_{n \geq 0} S_n \leq q^2$. As an application, we show that Syracuse sequences are "Branch sequences" with all the required conditions for convergence and therefore this confirms the Collatz conjecture.*

## 1. Introduction and preliminary definitions

Sequences of rational powers $(\xi(p/q)^n)_{n \geq 0}$, with $\xi$ a positive real number and $p, q$ two positive coprime integers, have been extensively studied [10] and have connection with many number theory and combinatorics famous problems such as Syracuse problem (or Collatz problem), Z-numbers problem (or Mahler problem) and Waring problem ... [10][11][12] [13] [14][15][16][17]. They are so inextricable that Paul Erdos stated: *"Mathematics may not be ready for such problems"* [21]. In this paper, we will focus our attention on a family of sequences which are perturbations of the sequences of rational powers under the form $((\xi + \Sigma_n)(p^n/q^{n+e_n}))_{n \geq 0}$, we call "Branch sequences". We establish that the influence of some particular perturbations $(e_n, \Sigma_n)_{n \geq 0}$ can be bounded in the associated rational powers sequence $(\xi(p/q)^n)_{n \geq 0}$ and this makes possible to demonstrate the convergence of Branch sequences under certain conditions. The proof uses only elementary arithmetic results and the study of asymptotical behavior of such sequences. As an application we show that Syracuse sequences can be formulated as Branch sequences which fulfill the required conditions for convergence.

**Definition 1.1 (Sequence of rational powers).** *We consider $(R_n)_{n \geq 0}$ a sequence defined by:*

$$R_n = \xi \left(\frac{p}{q}\right)^n$$

*With:*

- $\xi > 1$ *a positive real number and not a power of $q$.*

- $p > q > 1$ coprime integers, $p \wedge q = 1$. The euclidian division of $(p, q)$ gives: $p = \alpha q + \beta$, with $1 \leq \alpha$ and $1 \leq \beta < q - 1$.

- We suppose we have: $\alpha + \beta \leq q$
  In particular this imply:
  $$\begin{cases} 1 \leq \alpha, \beta \leq q - 1 \\ p \leq q^2 - \beta(q-1) \leq q^2 - 1 \\ \log_q\left(\frac{p}{q}\right) < 1 \end{cases}$$

□

We consider now a family of perturbations to the sequence $(R_n)_{n \geq 0}$ called "Branch sequences". The origin of the name is related to a visualization process of the sequence and will be explained at the end of the paper.

**Definition 1.2 (Branch sequence v1)** *We call Branch sequence associated to the sequence of rational powers $(R_n)_{n \geq 0}$ a sequence $(T_n)_{n \geq 0}$ defined by:*

- $T_0 = \frac{\xi}{q^{g_0}} > 1$ *with $g_0 \geq 0$ a positive integer.*
- *for $n \geq 0$ we take:*

$$T_{n+1} = \begin{cases} \frac{p(T_n + r_n)}{q} & \text{if } \left(\left\lfloor \frac{T_n}{q} \right\rfloor \equiv 0 \bmod q^2\right) \text{ or } \left(\left\lfloor \frac{T_n}{q} \right\rfloor \equiv (q^2 - 1) \bmod q^2\right) & (1.1.a) \\ \frac{T_n}{q} & \text{if } \left(\left\lfloor \frac{T_n}{q} \right\rfloor \not\equiv 0 \bmod q^2\right) \text{ and } \left(\left\lfloor \frac{T_n}{q} \right\rfloor \neq (q^2 - 1) \bmod q^2\right) & (1.1.b) \end{cases}$$

*With $(r_n)_{n \geq 0}$ real sequence chosen, for each n, such that:*

$$\lfloor T_n + r_n \rfloor = \lfloor T_n \rfloor \tag{1.2.a}$$

*Or equivalently,*

$$0 \leq \{T_n\} + r_n < 1 \Leftrightarrow -\{T_n\} \leq r_n < 1 - \{T_n\} \tag{1.2.b}$$

*$\{*\}$ and $\lfloor * \rfloor$ are respectively the fractional part and the integer part functions* □

**Remark 1.1**
a) For all $n \geq 0$ we have $T_n \geq 1$.
b) If $T_n$ verify the condition (1.1.b) then there exists $m > n$ such that $T_m$ verify (1.1.a).

Indeed: The first assertion can be proved easily by recurrence: we have $T_0 \geq 1$, if we suppose that it is still true until $n \geq 0$ then:
- If $T_n \geq q$ in both cases of *(Definition 1.2)* we have $T_{n+1} \geq 1$
- Else $1 \leq T_n < q$, then : $\left\lfloor \frac{T_n}{q} \right\rfloor = 0$, so equation (1.1.a) is verified, and we have:
$$T_{n+1} = \frac{p(T_n + r_n)}{q} \geq 1$$
(because $p > q$, $T_n \geq 1$ and $\lfloor T_n + r_n \rfloor = \lfloor T_n \rfloor$)

For the second assertion, if it is false then:
$T_{n+k}$ verifies the condition (1.1.b) for every $k \geq 0$:

$$\left(\left\lfloor\frac{T_n}{q^k}\right\rfloor \neq 0 \bmod q^2\right) \text{ and } \left(\left\lfloor\frac{T_n}{q^k}\right\rfloor \neq (q^2-1) \bmod q^2\right)$$

But for $k$ large enough we have $\left\lfloor\frac{T_n}{q^k}\right\rfloor = 0$ and this contradicts the above equation.

Another equivalent definition, more appropriate for our purpose and reduced to elements of $(T_n)_{n\geq 0}$ verifying (1.1.b), is as follows:

**Definition 1.3 (Branch sequence v2)** *We call Branch sequence $(S_n)_{n\geq 0}$ associated to the sequence of rational powers $(R_n)_{n\geq 0}$ a sequence defined by:*

- $S_0 = \frac{\xi}{q^{g_0}} > 1$  *with $g_0 \geq 0$ a positive integer.*

  *Without loss of generality we will choose $g_0 = 0$ (if not we replace $\xi$ by $\frac{\xi}{q^{g_0}}$), so we suppose $S_0$ verifying (1.1.a):*

  $$\left(\left\lfloor\frac{S_0}{q}\right\rfloor \equiv 0 \bmod q^2\right) \text{ or } \left(\left\lfloor\frac{S_0}{q}\right\rfloor \equiv (q^2-1) \bmod q^2\right)$$

- *For $n \geq 0$ we put:*
  $$S_{n+1} = \frac{p(S_n + r_n)}{q^{1+g_{n+1}}}$$
  $(r_n)_{n\geq 0}$ *is a real sequence chosen, for each $n$, such that (as in Definition 1.2)::*

  $$\lfloor S_n + r_n \rfloor = \lfloor S_n \rfloor \qquad (1.3.a)$$

  *Or equivalently*

  $$0 \leq \{S_n\} + r_n < 1 \Leftrightarrow -\{S_n\} \leq r_n < 1 - \{S_n\} \qquad (1.3.b)$$

  $(r_n)_{n\geq 0}$ *will be designated, as the perturbation sequence.*

  *We define:*
  $$g_{n+1} = \vartheta_q\left(\frac{p(S_n + r_n)}{q}\right)$$

  *With $\vartheta_q$ is a valuation function that count the smallest index for $S_{n+1}$ to verify the condition (1.1.a). This means that for $0 \leq k \leq g_{n+1} - 1$ and $S'_n = \frac{p(S_n + r_n)}{q}$ we have:*

  $$\left(\left\lfloor\frac{S'_n}{q^{k+1}}\right\rfloor \neq 0 \bmod q^2\right) \text{ and } \left(\left\lfloor\frac{S'_n}{q^{k+1}}\right\rfloor \neq (q^2-1) \bmod q^2\right) \qquad (1.3.c)$$

  *And (as $S_{n+1} = \frac{S'_n}{q^{g_{n+1}}}$):*
  $$\left(\left\lfloor\frac{S_{n+1}}{q}\right\rfloor \equiv 0 \bmod q^2\right) \text{ or } \left(\left\lfloor\frac{S_{n+1}}{q}\right\rfloor \equiv (q^2-1) \bmod q^2\right) \qquad (1.3.d)$$

  □

**Remark 1.2**
   a) The sequence $(S_n)_{n\geq 0}$ is well defined and $S_n \geq 1 \; \forall n \geq 0$.
   b) If for some $n \geq 0$ we have $1 \leq S_n < q$ then $\forall k \geq 0$ we have:
      $$1 \leq S_{n+k} < q \text{ and } g_{n+k} = 0 \text{ or } 1.$$
Indeed:

a) Because $(S_n)_{n\geq 0}$ is a subsequence of $(T_n)_{n\geq 0}$ and according to *(Remark 1.1)*, the sequence $(S_n)_{n\geq 0}$ is well defined and $S_n \geq 1\ \forall n \geq 0$.

b) If $1 \leq S'_n = \frac{p(S_n + r_n)}{q} < q$ then the condition (1.3.d) is verified for $S'_n$ and $g_n = 0$, so:
$$1 \leq S_{n+1} = \frac{p(S_n + r_n)}{q} < q$$
Otherwise we have $q \leq S'_n < q^2$ because:
$$1 \leq S_n < q,\ \lfloor S_n + r_n \rfloor = \lfloor S_n \rfloor \text{ and } p < q^2$$
In this case we have $\left\lfloor \frac{S'_n}{q^2} \right\rfloor = 0$ and $g_n = 1$ and we have:
$$1 \leq S_{n+1} = \frac{p(S_n + r_n)}{q^2} < q$$

The Branch sequence $(S_n)_{n\geq 0}$ can be considered as a perturbation of the sequence of rational powers $(\xi(p/q)^n)_{n\geq 0}$ by the introduction of $(e_n, \Sigma_n)_{n\geq 0}$ sequences as shown in the next lemma. This gives another definition of Branch sequence $(S_n)_{n\geq 0}$.

**Lemma 1.1** *From Definition (1.3) we have a direct expression of $S_n$ in terms of $(\xi, p, q, (r_k)_{k\leq n}, (g_k)_{k\leq n})$:*

$$\forall n \geq 0,\ S_n = \left(\frac{p^n}{q^{n+e_n}}\right)(\xi + \Sigma_n) \tag{1.4}$$

with

$$\begin{cases} \Sigma_0 = 0 \\ e_0 = g_0 = 0 \end{cases} \text{ and } \begin{cases} \Sigma_n = \sum_{j=0}^{n-1} r_j \frac{q^{j+e_j}}{p^j} \\ e_n = \sum_{j=0}^{n} g_j \end{cases}$$

**Proof.** By recurrence: It is evident for n=0:
$$S_0 = \left(\frac{p^0}{q^{0+e_0}}\right)(\xi + \Sigma_0)$$
Suppose it is true for n then for n+1 we have:
$$S_{n+1} = \frac{p(S_n + r_n)}{q^{1+g_{n+1}}}$$
$$= \frac{p}{q^{1+g_{n+1}}}\left(\left(\frac{p^n}{q^{n+e_n}}\right)(\xi + \Sigma_n) + r_n\right)$$
$$= \left(\frac{p^{n+1}}{q^{n+1+e_{n+1}}}\right)\left(\xi + \Sigma_n + r_n \frac{q^{n+e_n}}{p^n}\right)$$
$$= \left(\frac{p^{n+1}}{q^{n+1+e_{n+1}}}\right)(\xi + \Sigma_{n+1})$$
■

We note that $(e_n)_{n\geq 0}$ is a positive integer growing sequence, and if $(r_n)_{n\geq 0}$ is a positive sequence then $(\Sigma_n)_{n\geq 0}$ is also a real growing sequence.

**Definition 1.4 (Deterministic Branch sequences and positive controlled perturbations)**
*We say that a Branch sequence $(S_n)_{n\geq 0}$ is deterministic if:*

$$\forall n > m \geq 0, \lfloor S_n \rfloor = \lfloor S_m \rfloor \Rightarrow \begin{cases} S_{n+1} = S_{m+1} \\ r_{n+1} = r_{m+1} \end{cases} \tag{1.5.a}$$

*(we have then also: $g_{n+1} = g_{m+1}$)*

*We say that $(r_n)_{n\geq 0}$ is a positive controlled perturbation sequence if:*

$$\forall n \geq 0, 0 < r_n < 1 \text{ and } \limsup_{n\to\infty} \sqrt[n]{r_n} = 1 \qquad (1.5.b)$$

*We say that $(r_n)_{n\geq 0}$ is a perturbation sequence of Syracuse type if:*

$$\forall n \geq 0, r_n = \frac{c}{q^{g_n}} \qquad (1.5.c)$$

*(c > 0 is a positive constant)* □

The above conditions will be useful for establishing the convergence of Branch sequences. Indeed, with the deterministic condition the existence of $\lfloor S_n \rfloor = \lfloor S_m \rfloor$, for $n \neq m$, will imply that the sequence $(S_n)_{n\geq 0}$ will be periodical and with the controlled perturbation condition we can study the asymptotic behavior of the sequence $(\Sigma_n)_{n\geq 0}$ independently from the perturbation sequence $(r_n)_{n\geq 0}$.

Our main result is a convergence proof, in some sens and under some assumptions, of Branch sequences and the application of this result to Syracuse sequences.

**Theorem 1.1 (Main result)** *If $(S_n)_{n\geq 0}$ is a deterministic Branch sequences with perturbation of type Syracuse (Definition 1.3, Definition 1.4) then:*

$$\min_{n\geq 0} S_n \leq q^2. \qquad (1.6)$$

*In particular, Syracuse sequences can always be associated to Branch sequences with the required conditions and this confirms the Collatz conjecture.*

**Proof.** It is a direct consequence of *(Theorem 3.1) and (Theorem 4.1).*

■

In section 2, we explore some interesting properties of Branch sequences and we establish, due to the perturbations conditions, the independence of the Branch sequence integer part (after the second digit), at each step, from the perturbation choice in its interval. In particular this implies that, at each step n, the sum of perturbations has the same effect as the largest perturbation and this suggests that, in general, the perturbation is still confined in the "beginning" of the Branch sequence q-base expansion (Theorem 2.1). In section 3, we propose an asymptotic study of Branch sequence behavior and this leads to a fundamental theorem of convergence (Theorem 3.1) of Branch sequences under some particular conditions on perturbations. In section 4, as an application, we establish that Syracuse sequences can always be associated with Branch sequences with fulfilled convergence conditions and this establishes that the Syracuse conjecture is true (Theorem 4.1). Finally, in section 5 we present an interesting visualization of the Syracuse sequences as Cellular Automata (CA) in the gray code binary expansion. This allows us to see that the CA has a "tree" structure with coalescing "Branches" and explains, in a "picture", the convergence of Syracuse sequences. It also justifies the name of "Branch sequences".

## 2. Branch sequences Properties

**Definition 2.1** *Let $(R_n)_{n\geq 0}$ be a sequence of rational powers (Definition 1.1) and $(S_n)_{n\geq 0}$ a Branch sequence associated to $(R_n)_{n\geq 0}$ with the (Definition 1.3) notations and conditions:*

$$\begin{cases} \text{For } n \geq 0 : S_{n+1} = \frac{p(S_n + r_n)}{q^{1+g_n}} \quad Or \quad S_n = \left(\frac{p^n}{q^{n+e_n}}\right)(\xi + \Sigma_n) \\ \text{with: } \Sigma_n = \sum_{j=0}^{n-1} r_j \frac{q^{j+e_j}}{p^j}, \quad \Sigma_0 = 0, g_j = 0 \text{ and } e_n = \sum_{j=0}^{n} g_j \end{cases}$$

(2.1.a)

*We consider then the two sequences $(C_n)_{n\geq 0}$ and $(\Delta_n)_{n\geq 0}$ respectively defined by:*

$$\begin{cases} C_n = \frac{R_n}{q^{e_n}} = \frac{p^n \xi}{q^{n+e_n}} \\ \Delta_n = \frac{p^n \Sigma_n}{q^{n+e_n}} \end{cases}$$

(2.1.b)

*We have then:*

$$S_n = C_n + \Delta_n \tag{2.1.c}$$

$$\Delta_{n+1} = \frac{p(\Delta_n + r_n)}{q^{1+g_n}} \tag{2.1.d}$$

□

**Remark 2.1 (Some preliminary results about q-base expansion)**

Let's consider the q-base expansion of both $S_n$, $\frac{pS_n}{q}$, $(S_n + r_n)$ and $\frac{p(S_n + r_n)}{q}$:

$$\begin{cases} S_n = \sum_{j \in \mathbb{Z}} s_n^j q^j \\ \frac{pS_n}{q} = \sum_{j \in \mathbb{Z}} s'^j_n q^j \\ S_n + r_n = \sum_{j \in \mathbb{Z}} t_n^j q^j \\ \frac{p(S_n + r_n)}{q} = \sum_{j \in \mathbb{Z}} t'^j_n q^j \end{cases}$$

(2.3.a)

( with $0 \leq s_n^j, s'^j_n, t_n^j, t'^j_n \leq q - 1$ and we complete the q-expansion with zeros if necessary)

Then we have the next arithmetic elementary properties:

a) The arithmetic addition: $\alpha S_n + \beta \frac{S_n}{q} = \frac{pS_n}{q}$ gives the next relations between q-base expansion coefficients $(s_n^{j+1}, s_n^j \text{ and } s'^j_n)$ and the carry-digits $(\delta_n^j)_{j \in \mathbb{Z}}$

$$\beta s_n^{j+1} + \alpha s_n^j + \delta_n^j = s'^j_n + q \delta_n^{j+1} \tag{2.3.b}$$

$$\delta_n^j = \left\lfloor \frac{\beta s_n^j + \alpha s_n^{j-1} + \delta_n^{j-1}}{q} \right\rfloor \quad \text{and} \quad s'^j_n = \beta s_n^{j+1} + \alpha s_n^j + \delta_n^j \mod(q) \tag{2.3.c}$$

Relations (2.3.b) (2.3.c) for $n \geq 0$ are called the transition rules for $(S_n)_{n\geq 0}$.

b) By summing (2.3.b) on indexes $k \leq j$ we have:

$$\beta \sum_{k \leq j} s_n^k q^k + \alpha q \sum_{k \leq j} s_n^{k-1} q^{k-1} + q \sum_{k \leq j} \delta_n^{k-1} q^{k-1} = q \sum_{k \leq j} s_n'^{k-1} q^{k-1} + q \sum_{k \leq j} \delta_n^k q^k$$

Which gives: $\beta s_n^j + p \left\{\frac{S_n}{q^j}\right\} = q \left\{\frac{pS_n}{q^{j+1}}\right\} + q \delta_n^j$    i.e.:

$$\left\{\frac{pS_n}{q^{1+j}}\right\} = \frac{p}{q} \left\{\frac{S_n}{q^j}\right\} - \delta_n^j + \frac{\beta}{q} s_n^j \tag{2.3.d}$$

As $\left\{\frac{S_n}{q^j}\right\} < 1$ $and$ $\alpha + \beta \leq q$ we note that:

$$\delta_n^j = \left\lfloor \frac{\beta s_n^j + p\left\{\frac{S_n}{q^j}\right\}}{q} \right\rfloor \leq \left\lfloor \frac{\beta(q-1) + (\alpha q + \beta)\left\{\frac{S_n}{q^j}\right\}}{q} \right\rfloor \leq \left\lfloor \frac{(\alpha+\beta)q}{q} - \frac{\beta}{q}\left(1 - \left\{\frac{S_n}{q^j}\right\}\right) \right\rfloor < q \tag{2.3.e}$$

c) By summing (2.3.b) on indexes $k \geq j$ we have the relation:

$$\beta \sum_{k \geq j} s_n^k q^k + \alpha q \sum_{k \geq j} s_n^{k-1} q^{k-1} + q \sum_{k \geq j} \delta_n^{k-1} q^{k-1} = q \sum_{k \geq j} s_n'^{k-1} q^{k-1} + q \sum_{k \geq j} \delta_n^k q^k$$

Which gives: $-\beta s_n^j + p \left\lfloor \frac{S_n}{q^j} \right\rfloor + q \delta_n^j = q \left\lfloor \frac{pS_n}{q^{1+j}} \right\rfloor$ i.e. :

$$\left\lfloor \frac{pS_n}{q^{1+j}} \right\rfloor = \left\lfloor \frac{pS_n}{q^{1+j}} \right\rfloor = \frac{p}{q} \left\lfloor \frac{S_n}{q^j} \right\rfloor + \delta_n^j - \frac{\beta}{q} s_n^j \tag{2.3.f}$$

d) We have equivalently the same transition rules with $(S_n + r_n)$ and $\frac{p(S_n + r_n)}{q}$

$$\beta t_n^{j+1} + \alpha t_n^j + \sigma_n^j = t_n'^j + q \sigma_n^{j+1} \tag{2.3.g}$$

$$\sigma_n^j = \left\lfloor \frac{\beta t_n^j + \alpha t_n^{j-1} + \sigma_n^{j-1}}{q} \right\rfloor \quad \text{and} \quad t_n'^j = \beta t_n^{j+1} + \alpha t_n^j + \sigma_n^j \bmod(q) \tag{2.3.h}$$

Which gives:

$$\left\{\frac{p(S_n + r_n)}{q^{1+j}}\right\} = \frac{p}{q} \left\{\frac{(S_n + r_n)}{q^j}\right\} - \sigma_n^j + \frac{\beta}{q} t_n^j \tag{2.3.i}$$

$$\left\lfloor \frac{p(S_n + r_n)}{q^{1+j}} \right\rfloor = \frac{p}{q} \left\lfloor \frac{(S_n + r_n)}{q^j} \right\rfloor + \sigma_n^j - \frac{\beta}{q} t_n^j \tag{2.3.j}$$

$\left(\sigma_n^j\right)_{j \in \mathbb{Z}}$ is the carry digits sequence associated to the arithmetic addition:

$$\alpha(S_n + r_n) + \beta \frac{(S_n + r_n)}{q} = \frac{p(S_n + r_n)}{q}$$

                                                                                              □

Now we will establish an important result for the rest of the paper: the independence of the carry digits $(\delta_n^k)_{k \geq 2}$ (in the nth transition rule of $(S_n)_{n \geq 0}$) from the choice of $r_n \in [-\{S_n\}, 1 - \{S_n\}[$.

**Lemma 2.1 (Independence of carry digits for $(S_n)_{n \geq 0}$)** *For $n \geq 0$ and $k \geq 2$ we have:*

$$\left\lfloor \frac{p(S_n + r_n)}{q^{1+k}} \right\rfloor = \left\lfloor \frac{pS_n}{q^{1+k}} \right\rfloor = \frac{p}{q} \left\lfloor \frac{S_n}{q^k} \right\rfloor + \delta_n^k - \frac{\beta s_n^k}{q} \tag{2.4.a}$$

*In other words we have for $j \geq k$: $t'^j_n = s'^j_n$ and $\sigma_n^j = \delta_n^j$ are independent of the choice of $r_n \in [-\{S_n\}, 1 - \{S_n\}[$.*

**Proof.** From (1.3.a) and (2.3.a) we have

$$\lfloor S_n + r_n \rfloor = \lfloor S_n \rfloor \xRightarrow{\text{for } j \geq 0} t_n^j = s_n^j \tag{2.4.b}$$

From (2.3.g) and (2.3.h) as $t_n^j = s_n^j$ for $j \geq 0$, we have:

$$\begin{cases} t'^k_n + q\sigma_n^{k+1} = \beta s_n^{k+1} + \alpha s_n^k + \sigma_n^k \\ \sigma_n^{k+1} = \left\lfloor \frac{\beta s_n^{k+1} + \alpha s_n^k + \sigma_n^k}{q} \right\rfloor \end{cases}$$

When $r_n$ is varying in: $[-\{S_n\}, 1 - \{S_n\}[$ we see that:

- $\sigma_n^{k+1}$ depends uniquely on $\sigma_n^k$ (as $s_n^{k+1}$ and $s_n^k$ still invariant)
- $t'^k_n$ is determined from $t'^k_n + q\sigma_n^{k+1}$ by uniqueness of q-base expansion

So it is sufficient to prove that $\sigma_n^2$ is independent of the choice of $r_n$ to have the same result for $\sigma_n^k$ for $k \geq 2$.

Now we have: $\quad \sigma_n^2 = \left\lfloor \frac{\beta s_n^2 + \alpha s_n^1 + \sigma_n^1}{q} \right\rfloor$

And: $\quad \left( \left\lfloor \frac{S_n}{q} \right\rfloor \equiv 0 \bmod(q^2) \right)$ or $\left( \left\lfloor \frac{S_n}{q} \right\rfloor \equiv (q^2 - 1) \bmod(q^2) \right)$

i.e.: $\quad s_n^1 = s_n^2 = 0$ or $s_n^1 = s_n^2 = q - 1$

So we have:

$$\begin{cases} \sigma_n^2 = \left\lfloor \frac{\sigma_n^1}{q} \right\rfloor = 0 & \text{if } \left( \left\lfloor \frac{S_n}{q} \right\rfloor \equiv 0 \bmod(q^2) \right) \\ \sigma_n^2 = \left\lfloor (q-1) + \frac{\sigma_n^1}{q} \right\rfloor = q - 1 & \text{if } \left( \left\lfloor \frac{S_n}{q} \right\rfloor \equiv 0 \bmod(q^2) \right) \end{cases}$$

Conclusion: in both cases $\sigma_n^2$ is independent of the choice of $r_n \in [-\{S_n\}, 1 - \{S_n\}[$

If we choose $r_n = 0$ we have then $\sigma_n^2 = \delta_n^2$ and this imply:

$$\sigma_n^k = \delta_n^k \text{ for } k \geq 2 \tag{2.4.c}$$

Now we have *for* $k \geq 2$:

$$\left\lfloor \frac{p(S_n+r_n)}{q^{1+k}} \right\rfloor = \frac{p}{q} \left\lfloor \frac{(S_n+r_n)}{q^k} \right\rfloor + \sigma_n^k - \frac{\beta}{q} t_n^k \text{ from (2.3.j)}$$

$$= \frac{p}{q} \left\lfloor \frac{(S_n+r_n)}{q^k} \right\rfloor + \delta_n^k - \frac{\beta}{q} s_n^k \text{ from (2.4.b) and (2.4.c)}$$

$$= \frac{p}{q} \left\lfloor \frac{S_n}{q^k} \right\rfloor + \delta_n^k - \frac{\beta}{q} s_n^k \text{ from (1.3.a)}$$

$$= \left\lfloor \frac{pS_n}{q^{1+k}} \right\rfloor \text{ from (2.3.f)}$$

■

The next lemma establishes an elementary arithmetic property of the integer part function $\lfloor * \rfloor$. As elementary as this result sounds, to the best of my knowledge, there is no record of it elsewhere. It will be useful to transfer the result of condition (1.3.1) and (Lemma 2.1) from the Branch sequence $(S_n)_{n \geq 0}$ to the sequence $(\Delta_n)_{n \geq 0}$.

Let's consider the q-base expansion of the positive numbers A , B and B':

$$\begin{cases} A = \sum_{j \in \mathbb{Z}} a_j q^j \\ B = \sum_{j \in \mathbb{Z}} b_j q^j \\ B' = \sum_{j \in \mathbb{Z}} b'_j q^j \end{cases} \quad (2.5.a)$$

( with $0 \leq a_n^j, b_n^j, b'_n^j \leq q - 1$ and we complete the q-expansion with zeros if necessary)

**Lemma 2.2** *For the positive numbers **A** , **B** and **B'** we have the next arithmetic elementary property:*

$$\text{If } \begin{cases} B < B' \\ \lfloor A + B \rfloor = \lfloor A + B' \rfloor \end{cases} \text{ then } \lfloor B \rfloor = \lfloor B' \rfloor \quad (2.5.b)$$

**Proof.** $\lfloor A + B \rfloor = \lfloor A + B' \rfloor$ is equivalent to:

$$\forall k \geq 0, \left\lfloor \frac{A}{q^k} + \frac{B}{q^k} \right\rfloor = \left\lfloor \frac{A}{q^k} + \frac{B'}{q^k} \right\rfloor$$

So:

$$\lfloor A + B \rfloor = \lfloor A + B' \rfloor \Leftrightarrow \lfloor B \rfloor + \lfloor \{A\} + \{B\} \rfloor = \lfloor B' \rfloor + \lfloor \{A\} + \{B'\} \rfloor$$

$$\Leftrightarrow \forall k \geq 0, \left\lfloor \frac{B}{q^k} \right\rfloor + \left\lfloor \left\{\frac{A}{q^k}\right\} + \left\{\frac{B}{q^k}\right\} \right\rfloor = \left\lfloor \frac{B'}{q^k} \right\rfloor + \left\lfloor \left\{\frac{A}{q^k}\right\} + \left\{\frac{B'}{q^k}\right\} \right\rfloor \quad (2.5.c)$$

As $B < B'$ , let $k$ be the least index such that: $b_k < b'_k$ then:

$$\begin{cases} \forall k' > k: \left\lfloor \frac{B}{q^{k'}} \right\rfloor = \left\lfloor \frac{B'}{q^{k'}} \right\rfloor (\Leftrightarrow \forall k' > k \ b_{k'} = b'_{k'}) \\ \left\lfloor \frac{B}{q^k} \right\rfloor < \left\lfloor \frac{B'}{q^k} \right\rfloor \end{cases}$$

If $k \geq 0$ then by (2.5.c) we have then:

$$\left\lfloor \frac{B}{q^k} \right\rfloor < \left\lfloor \frac{B'}{q^k} \right\rfloor \Rightarrow \left\lfloor \left\{\frac{A}{q^k}\right\} + \left\{\frac{B}{q^k}\right\} \right\rfloor > \left\lfloor \left\{\frac{A}{q^k}\right\} + \left\{\frac{B'}{q^k}\right\} \right\rfloor$$

But since we have:

$$0 \leq \left\lfloor \left\{\frac{A}{q^k}\right\} + \left\{\frac{B'}{q^k}\right\} \right\rfloor < \left\lfloor \left\{\frac{A}{q^k}\right\} + \left\{\frac{B}{q^k}\right\} \right\rfloor \leq 1$$

Then necessarily:

$$\begin{cases} \left\lfloor \left\{\frac{A}{q^k}\right\} + \left\{\frac{B}{q^k}\right\} \right\rfloor = 1 \\ \left\lfloor \left\{\frac{A}{q^k}\right\} + \left\{\frac{B'}{q^k}\right\} \right\rfloor = 0 \end{cases}$$

If we consider the q-base expansion of $A + B$ and $A + B'$:

$$\begin{cases} A + B = \sum_{j \in \mathbb{Z}} c_j q^j \\ A + B' = \sum_{j \in \mathbb{Z}} c'_j q^j \end{cases}$$

Then:

$$\begin{cases} c_{k+1} = a_{k+1} + b_{k+1} + \left\lfloor \left\{\frac{A}{q^k}\right\} + \left\{\frac{B}{q^k}\right\} \right\rfloor \bmod(q) \\ c'_{k+1} = a_{k+1} + b'_{k+1} + \left\lfloor \left\{\frac{A}{q^k}\right\} + \left\{\frac{B'}{q^k}\right\} \right\rfloor \bmod(q) \end{cases}$$

So as: $b_{k+1} = b'_{k+1}$ we have:

$$c_{k+1} = c'_{k+1} + 1 \bmod(q) \Rightarrow c_{k+1} \neq c'_{k+1}$$

Which contradict, as $k \geq 0$, the fact that $\lfloor A + B \rfloor = \lfloor A + B' \rfloor$

So we have necessarily $k < 0$ and then:    $\lfloor B \rfloor = \lfloor B' \rfloor$

■

Due to this elementary arithmetic result we can retrieve for $(\Delta_n)_{n \geq 0}$ the same result about the carry digits independence obtained for $(S_n)_{n \geq 0}$.

**Corollary 2.1 (Independence of carry digits for $(\Delta_n)_{n \geq 0}$ )** *If the perturbation sequence $(r_n)_{n \geq 0}$ is positive (the sequence $(\Delta_n)_{n \geq 0}$ is then also positive) then for $n \geq 0$ and $k \geq 2$ we have:*

$$\begin{cases} \lfloor \Delta_n + r_n \rfloor = \lfloor \Delta_n \rfloor \\ \left\lfloor \frac{p(\Delta_n + r_n)}{q^{1+k}} \right\rfloor = \left\lfloor \frac{p\Delta_n}{q^{1+k}} \right\rfloor \end{cases}$$

**Proof.** It is a direct application of Lemma 2.1 and Lemma 2.2 :

$$\lfloor S_n + r_n \rfloor = \lfloor S_n \rfloor \Leftrightarrow \lfloor C_n + (\Delta_n + r_n) \rfloor = \lfloor C_n + \Delta_n \rfloor$$
$$\Rightarrow \lfloor \Delta_n + r_n \rfloor = \lfloor \Delta_n \rfloor$$

and

$$\left\lfloor \frac{p(S_n+r_n)}{q^{1+k}} \right\rfloor = \left\lfloor \frac{pS_n}{q^{1+k}} \right\rfloor \Leftrightarrow \left\lfloor \frac{pC_n}{q^{1+k}} + \frac{p(\Delta_n+r_n)}{q^{1+k}} \right\rfloor = \left\lfloor \frac{pC_n}{q^{1+k}} + \frac{p\Delta_n}{q^{1+k}} \right\rfloor$$
$$\Rightarrow \left\lfloor \frac{p(\Delta_n+r_n)}{q^{1+k}} \right\rfloor = \left\lfloor \frac{p\Delta_n}{q^{1+k}} \right\rfloor$$

∎

Now we will define three new sequences to help us study the convergence of the branch sequences. In those three sequences we will replace the sum of perturbations $\sum_{j=0}^{n-1} r_j \frac{q^{j+e_j}}{p^j}$ by the largest perturbation $max_{0 \leq j \leq n-1}\left(r_j \frac{q^{j+e_j}}{p^j}\right)$.

**Definition 2.2** *We define the sequences* $(\omega_n)_{n \geq 0}$, $(\Omega_n)_{n \geq 0}$ *and* $(Z_n)_{n \geq 0}$ *(with the Definition 2.1 notations) respectively by,* $\forall n \geq 0$:

$$\omega_n = max_{0 \leq j \leq n-1}\left(r_j \frac{q^{j+e_j}}{p^j}\right)$$

$$\Omega_n = \frac{p^n}{q^{n+e_n}} max_{0 \leq j \leq n-1}\left(r_j \frac{q^{j+e_j}}{p^j}\right) = \frac{p^n}{q^{n+e_n}} \omega_n$$

$$Z_n = \frac{p^n}{q^{n+e_n}}(\xi + \omega_n)$$

The three sequences have respectively the same form as $(\Sigma_n)_{n \geq 0}$, $(\Delta_n)_{n \geq 0}$ and $(S_n)_{n \geq 0}$ with the difference of replacing the sum of perturbations $\sum_{j=0}^{n-1} r_j \frac{q^{j+e_j}}{p^j}$ by the largest perturbation $max_{0 \leq j \leq n-1}\left(r_j \frac{q^{j+e_j}}{p^j}\right)$.                     □

**Lemma 2.3** *If the perturbation sequence* $(r_n)_{n \geq 0}$ *is positive, then sequence* $(\Omega_n)_{n \geq 0}$ *(Definition 2.2) verify for* $n \geq 0$ *and* $k \geq 2$:

$$\left\lfloor \frac{\Omega_{n+1}}{q^k} \right\rfloor = \left\lfloor \frac{p\Omega_n}{q^{1+g_{n+1}+k}} \right\rfloor$$

**Proof.** By definition of $(\Omega_n)_{n \geq 0}$ we have:

$$\Omega_{n+1} = \frac{p^{n+1}}{q^{n+1+e_{n+1}}} max_{0 \leq j \leq n}\left(r_j \frac{q^{j+e_j}}{p^j}\right)$$
$$= \frac{p}{q^{1+g_{n+1}}} max_{0 \leq j \leq n}\left(r_j \frac{p^{n-j}}{q^{n-j+e_n-e_j}}\right)$$
$$= \frac{p}{q^{1+g_{n+1}}} max\left(r_n, \frac{p^n}{q^{n+e_n}} max_{0 \leq j \leq n-1}\left(r_j \frac{q^{j+e_j}}{p^j}\right)\right)$$
$$= \frac{p}{q^{1+g_{n+1}}} max(r_n, \Omega_n)$$

So we have two cases:

· If $r_n \geq \Omega_n$ then:

$$\Omega_{n+1} = \frac{p}{q^{1+g_{n+1}}} r_n$$

we have $r_n < 1$ ( Definition 1.3 ) an this imply also $\Omega_n < 1$:

And as $p < q^2$ (Definition 1.1) we have for $k \geq 2$:

$$\frac{p}{q^{1+g_{n+1}+k}} < 1 \text{ and then } \begin{cases} \frac{p}{q^{1+g_{n+1}+k}} \Omega_n < 1 \\ \frac{\Omega_{n+1}}{q^k} = \frac{p}{q^{1+g_{n+1}+k}} r_n < 1 \end{cases}$$

So in the case $r_n \geq \Omega_n$ we have :

$$\left[\frac{\Omega_{n+1}}{q^k}\right] = \left[\frac{p\Omega_n}{q^{1+g_{n+1}+k}}\right] = 0 \tag{2.6}$$

- If $r_n < \Omega_n$ then:

$$\Omega_{n+1} = \frac{p}{q^{1+g_{n+1}}} \Omega_n$$

And of course :

$$\left[\frac{\Omega_{n+1}}{q^k}\right] = \left[\frac{p\Omega_n}{q^{1+g_{n+1}+k}}\right]$$

∎

The next theorem establishes a fundamental link between the sequences $(\Delta_n)_{n\geq 0}$ and $(\Omega_n)_{n\geq 0}$ which will open the way to the convergence study of the branch sequences $(S_n)_{n\geq 0}$. It states that the perturbation, represented by the sequence $(\Delta_n)_{n\geq 0}$, is always "dominated" by the largest perturbation sequence, represented by $(\Omega_n)_{n\geq 0}$. This will suggest that the perturbation sequences, provoked by Branch sequences, remains in general confined in the beginning digits of the rational powers sequence q-base expansion as, each time the largest perturbation is reached, we have $\Omega_n \leq q^2$ and then $\Delta_n \leq q^2$.

### Theorem 2.1 (Largest perturbation "domination")
*If the perturbation sequence $(r_n)_{n\geq 0}$ is positive then for $n \geq 1$ and $k \geq 2$ we have:*

$$\left[\frac{\Delta_n}{q^k}\right] = \left[\frac{\Omega_n}{q^k}\right] \tag{2.7}$$

**Proof.** Let's consider the q-base expansion, *for $n \geq 0$*, of both $\Delta_n$, $\frac{p\Delta_n}{q}$, $(\Delta_n + r_n)$ and $\Delta'_n = \frac{p(\Delta_n + r_n)}{q}$:

$$\begin{cases} \Delta_n = \sum_{j\in Z} d_n^j q^j \\ \frac{p\Delta_n}{q} = \sum_{j\in Z} d'^j_n q^j \\ \Delta_n + r_n = \sum_{j\in Z} z_n^j q^j \\ \frac{p(\Delta_n + r_n)}{q} = \sum_{j\in Z} z'^j_n q^j \end{cases}$$

And the transition rules (Remark 2.1) of $(\Delta_n)_{n\geq 0}$ and $(\Delta_n + r_n)_{n\geq 0}$ *for $j \in Z$:*

$$\begin{cases} \beta d_n^{j+1} + \alpha d_n^j + \tau_n^j = d'_n^j + q\tau_n^{j+1} \\ \beta z_n^{j+1} + \alpha z_n^j + \tau'_n^j = z'_n^j + q\tau'_n^{j+1} \end{cases}$$

( $(\tau_n^j)_{j\in Z}$ and $(\tau'_n^j)_{j\in Z}$ are the carry digits (Remark 2.1) )

We note that by (Definition 2.1) and (Corollary 2.1) we have:

- $\forall j \in Z, d'_n^j = d_{n+1}^{j-g_{n+1}}$ as $\Delta_{n+1} = \frac{p(\Delta_n + r_n)}{q^{1+g_n}}$ (2.8.a)
- $\forall j \geq 0, d_n^j = z_n^j$ as $\lfloor \Delta_n + r_n \rfloor = \lfloor \Delta_n \rfloor$ (2.8.b)
- $\forall k \geq 2, d'^k_n = z'^k_n$ as $\left\lfloor \frac{p(\Delta_n+r_n)}{q^{1+k}} \right\rfloor = \left\lfloor \frac{p\Delta_n}{q^{1+k}} \right\rfloor$ (2.8.c)
- $\forall k \geq 2, \tau_n^k = \tau'^k_n$ as $\begin{cases} \tau_n^k = d'^k_n - \beta d_n^{k+1} - \alpha d_n^k \mod(q) \\ \tau'^k_n = z'^k_n - \beta z_n^{k+1} - \alpha z_n^k \mod(q) \end{cases}$ (2.8.d)

  It means, in particular, as noted in Lemma 2.1 and Corollary 2.1, that $\tau'^k_n$ is independent of the choice of $r_n \in [-\{S_n\}, 1 - \{S_n\}[$.

We will proceed by recurrence, the relation (2.7) is true for n=1:

$$\left\lfloor \frac{\Delta_1}{q^k} \right\rfloor = \left\lfloor \frac{p(\Delta_0 + r_0)}{q^{1+g_1+k}} \right\rfloor$$
$$= \left\lfloor \frac{p}{q^{1+e_1+k}} r_0 \right\rfloor \quad (as\ \Delta_0 = 0\ and\ e_1 = g_1)$$
$$= \left\lfloor \frac{\Omega_1}{q^k} \right\rfloor$$

(we notice that $\left\lfloor \frac{\Delta_1}{q^k} \right\rfloor = \left\lfloor \frac{\Omega_1}{q^k} \right\rfloor = 0$ as $\frac{p}{q^{1+e_1+k}} r_0 < 1$ (Lemma 2.3)) (2.8.e)

Suppose it still true until the index $n \geq 1$ :

$$\left\lfloor \frac{\Delta_m}{q^k} \right\rfloor = \left\lfloor \frac{\Omega_m}{q^k} \right\rfloor \ for: 1 \leq m \leq n$$

Let's verify that it is true for the index $n + 1$: For k ≥ 2 and from (Corollary 2.1) and (Remark 2.1) and notes above, we have:

$$\left\lfloor \frac{\Delta_{n+1}}{q^k} \right\rfloor = \left\lfloor \frac{p(\Delta_n+r_n)}{q^{1+g_{n+1}+k}} \right\rfloor$$
$$= \left\lfloor \frac{p\Delta_n}{q^{1+g_{n+1}+k}} \right\rfloor \qquad (2.8.f)$$
$$= \frac{p}{q}\left\lfloor \frac{\Delta_n}{q^{g_{n+1}+k}} \right\rfloor + \tau_n^{g_{n+1}+k} - \frac{\beta d_n^{g_{n+1}+k}}{q}$$

Now by recurrence hypothesis we have:

$$\left\lfloor \frac{\Delta_n}{q^k} \right\rfloor = \left\lfloor \frac{\Omega_n}{q^k} \right\rfloor \Rightarrow \left\lfloor \frac{\Delta_n}{q^{g_{n+1}+k}} \right\rfloor = \left\lfloor \frac{\Omega_n}{q^{g_{n+1}+k}} \right\rfloor$$

Hence from (2.8.f):

$$\left\lfloor \frac{\Delta_{n+1}}{q^k} \right\rfloor = \frac{p}{q}\left\lfloor \frac{\Omega_n}{q^{g_{n+1}+k}} \right\rfloor + \tau_n^{g_{n+1}+k} - \frac{\beta d_n^{g_{n+1}+k}}{q}$$

Now let's consider the q-base expansion of $\Omega_n$:

$$\Omega_n = \sum_{j \in \mathbb{Z}} d''^j_n q^j$$

And its associated transition rules (Remark 2.1), $for\ j \in \mathbb{Z}$:

$$\beta d''^{j+1}_n + \alpha d''^j_n + \tau''^j_n = d''^{j-g_{n+1}}_{n+1} + q\tau''^{j+1}_n$$

( $\left(\tau''^j_n\right)_{j \in \mathbb{Z}}$ are the carry digits of the transition rules of $(\Omega_n)_{n \geq 0}$ )

($d''^{j-g_{n+1}}_{n+1}$ is a coefficient of the q-base expansion of $\Omega_{n+1}$)

This gives the relation (by (2.3.f)):

$$\lfloor \Omega_{n+1} \rfloor = \frac{p}{q}\left\lfloor \frac{\Omega_n}{q^{g_{n+1}+k}} \right\rfloor + \tau''^{g_{n+1}+k}_n - \frac{\beta d''^{g_{n+1}+k}_n}{q}$$

Again from (2.8.f) we have :

$$\left\lfloor \frac{\Delta_{n+1}}{q^k} \right\rfloor = \left\lfloor \frac{\Omega_n}{q^{g_{n+1}+k}} \right\rfloor - \tau''^{g_{n+1}+k}_n - \frac{\beta d''^{g_{n+1}+k}_n}{q} + \tau_n^{g_{n+1}+k} - \frac{\beta d_n^{g_{n+1}+k}}{q}$$

And, by recurrence hypothesis $\left\lfloor \frac{\Delta_n}{q^k} \right\rfloor = \left\lfloor \frac{\Omega_n}{q^k} \right\rfloor$, we have:

$$d''^{g_{n+1}+k}_n = d_n^{g_{n+1}+k}$$

So we have finally :

$$\left\lfloor \frac{\Delta_{n+1}}{q^k} \right\rfloor = \left\lfloor \frac{\Omega_{n+1}}{q^k} \right\rfloor + \left(\tau_n^{g_{n+1}+k} - \tau''^{g_{n+1}+k}_n\right)$$

Now the important fact is that from (2.8.d) and (Corollary 2.1) we have $\tau_n^{g_{n+1}+k}$ is independent of the choice of:

$$r_n \in [-\{S_n\}, 1 - \{S_n\}[$$

So if we choose in particular:

$$r_n = \left\{\frac{\Omega_n}{q^{g_{n+1}+k}}\right\} - \{S_n\}$$

We have :

$$-\{S_n\} \leq r_n = \left\{\frac{\Omega_n}{q^{g_{n+1}+k}}\right\} - \{S_n\} < 1 - \{S_n\}$$

and in this particular case of $r_n$ we have:

$$\tau_n^{g_{n+1}+k} = \tau''_n{}^{g_{n+1}+k}$$

So finally we can conclude that:

$$\left[\frac{\Delta_{n+1}}{q^k}\right] = \left[\frac{\Omega_{n+1}}{q^k}\right] \qquad \blacksquare$$

**Corollary 2.2** *If the perturbation sequence* $(r_n)_{n\geq 0}$ *is positive, then for* $n \geq 0$ :

$$\Omega_n < \Delta_n < \Omega_n + q^2 \qquad (2.9.a)$$

$$Z_n < S_n < Z_n + q^2 \qquad (2.9.b)$$

*(Or equivalently:* $C_n + \Omega_n \leq S_n < C_n + \Omega_n + q^2$*)*

$$\omega_n < \Sigma_n < \omega_n + \frac{q^{n+e_n}}{p^n}q^2 \qquad (2.9.c)$$

**Proof.** It is a direct application of Theorem 2.1.

As , (Definition 2.1, definition 2.2) , for n ≥ 1 :

$$S_n = C_n + \Delta_n \quad \text{with: } C_n = \frac{p^n}{q^{n+e_n}}\xi \text{ and } \Delta_n = \frac{p^n}{q^{n+e_n}}\Sigma_n$$

$$\Omega_n = \frac{p^n}{q^{n+e_n}}\omega_n \text{ and } Z_n = \frac{p^n}{q^{n+e_n}}(\xi + \omega_n)$$

we have $\Omega_n < \Delta_n$ which imply also:

$$\omega_n < \Sigma_n$$

$$Z_n < S_n$$

And from (Theorem 2.1 ) :

$$\left[\frac{\Delta_n}{q^2}\right] = \left[\frac{\Omega_n}{q^2}\right] \Rightarrow \Delta_n < \Omega_n + q^2$$

which imply also :

$$\Sigma_n < \omega_n + \frac{q^{n+e_n}}{p^n}q^2$$

$$S_n < Z_n + q^2$$

And we get directly  (2.9.a), (2.9.b)  and (2.9.c)

$$\blacksquare$$

## 3. Branch sequences convergence study

In this section we will analyze the asymptotic behavior of the sequence $(S_n)_{n\geq 0}$ and establish some convergence results. In the remaining part of the paper we will suppose that the branch sequence $(S_n)_{n\geq 0}$ is deterministic and have positive controlled perturbation (Definition 1.3, Definition 1.4).

**Lemma 3.1** *We have three possibilities for the sequence $(\Sigma_n)_{n\geq 0}$:*

- *if $limsup_{n\to\infty} \frac{e_n}{n} > log_q\left(\frac{p}{q}\right)$ then $(\Sigma_n)_{n\geq 0}$ is divergent and $lim_{n\to\infty} \Sigma_n = +\infty$*
- *if $limsup_{n\to\infty} \frac{e_n}{n} > log_q\left(\frac{p}{q}\right)$ then $(\Sigma_n)_{n\geq 0}$ is convergent and $lim_{n\to\infty} \Sigma_n = \Sigma < +\infty$*
- *if $limsup_{n\to\infty} \frac{e_n}{n} = log_q\left(\frac{p}{q}\right)$ then the situation is inconclusive.*

**Proof.** Firstly we note that $\Sigma_n = \sum_{j=0}^{n-1} r_j \frac{q^{j+e_j}}{p^j}$ is a partial sum of the power series:

$$\Sigma(z) = \sum_{n\geq 0} r_n q^{e_n} z^n \text{ for the value } z = \frac{p}{q}$$

The convergence radius R of the power series is given by the Cauchy-Hadamard theorem:

$$\frac{1}{R} = limsup_{n\to\infty} \sqrt[n]{r_n q^{e_n}}$$

As $limsup_{n\to\infty} \sqrt[n]{r_n} = 1$ (a controlled perturbation sequence) and $limsup_{n\to\infty} \sqrt[n]{q^{e_n}} > 0$ we have:

$$limsup_{n\to\infty} \sqrt[n]{r_n q^{e_n}} = limsup_{n\to\infty} \sqrt[n]{r_n} * limsup_{n\to\infty} \sqrt[n]{q^{e_n}}$$

So

$$\frac{1}{R} = limsup_{n\to\infty} \sqrt[n]{q^{e_n}}$$

Or by passing to logarithm:

$$log_q\left(\frac{1}{R}\right) = limsup_{n\to\infty} \frac{e_n}{n}$$

The lemma follows directly from the Cauchy-Hadamard criterion on power series convergence.

∎

**Lemma 3.2** *We have two exclusive cases for the sequence $(S_n)_{n\geq 0}$ :*
   a) *$(S_n)_{n\geq 0}$ is a bounded sequence and in this case we have:*
   - *$(S_n)_{n\geq 0}$ is periodic, with period m, after a certain order $m_0$ i.e.:*
     *$\exists m_0 \geq 0 \text{ and } m > 0, \forall k \geq 0 \text{ and } 0 \leq l < m, S_{m_0+km+l} = S_{m_0+l}$*
   - *$lim_{n\to\infty} \left(\frac{e_n}{n}\right) = \frac{e_m}{m}$*
   - *$lim_{n\to\infty} \frac{e_n}{n} > log_q\left(\frac{p}{q}\right)$*

- $\lim\limits_{n\to\infty} \Sigma_n = +\infty$

*Without loss of generality in this case we will always suppose that $m_0 = 0$ i.e. $(S_n)_{n\geq 0}$ is fully periodic (by for example beginning the sequence at $S_{m_0}$).*

b) $(S_n)_{n\geq 0}$ is an unbounded sequence and in this case:
- $\lim\limits_{n\to\infty} S_n = +\infty$
- $\limsup_{n\to\infty} \left(\frac{e_n}{n}\right) \leq \log_q\left(\frac{p}{q}\right)$

**Proof.**

a) If $(S_n)_{n\geq 0}$ is bounded then $(\lfloor S_n \rfloor)_{n\geq 0}$ is also bounded and by "Pigeonhole principle" we have :
$$\exists m_0 > m_1 \geq 0, \lfloor S_{m_0} \rfloor = \lfloor S_{m_1} \rfloor$$
Then as $(S_n)_{n\geq 0}$ is deterministic (Definition 1.4) we have:
$$g_{m_0+1} = g_{m_1+1} \text{ and } S_{m_0+1} = S_{m_1+1} \text{ and again } \lfloor S_{m_0+1} \rfloor = \lfloor S_{m_1+1} \rfloor$$
So for $k \geq 0$ and $0 \leq l < m$ with $m = m_1 - m_0$ :
$$S_{m_0+1+km+l} = S_{m_0+1+l} \text{ and } g_{m_0+km+l+1} = g_{m_0+1+l+1}$$

This means that the sequences $(S_n)_{n\geq 0}$ and $(g_{n+1})_{n\geq 0}$ are periodic from the index $n = m_0 + 1$. As we are interested in the asymptotic behavior, we can replace, without loss of generality, $S_0$ by $S_{m_0+1}$. Then $(S_n)_{n\geq 0}$ is fully periodic of period $m \geq 1$. This implies that $\forall k \geq 0$ and $0 \leq l < m$:

$$\begin{cases} S_{km+l} = S_l \\ e_{km+l} = k * e_m + e_l \end{cases}$$

Now for $n = km + l$ with $k \geq 0$ and $0 \leq l < m$ we have
$$\frac{e_n}{n} = \frac{e_{km+l}}{km+l} = \frac{k*e_m + e_l}{km+l} = \frac{e_m}{m+\frac{l}{k}} + \frac{e_l}{km+l} \xrightarrow[k\to+\infty]{} \frac{e_m}{m}$$

So we deduce:
$$\lim_{n\to\infty} \frac{e_n}{n} = \frac{e_m}{m}$$

By passing to q-base logarithm we have:
$$S_n = \left(\frac{p^n}{q^{n+e_n}}\right)(S_0 + \Sigma_n) \Leftrightarrow \frac{e_n}{n} = \log_q\left(\frac{p}{q}\right) + \frac{\log_q\left(1+\frac{\Sigma_n}{S_0}\right)}{n} + \frac{\log_q(S_0) - \log_q(S_n)}{n}$$

So we have:
$$\frac{e_m}{m} = \log_q\left(\frac{p}{q}\right) + \frac{\log_q\left(1+\frac{\Sigma_m}{S_0}\right)}{m} + \frac{\log_q(S_0) - \log_q(S_m)}{m} > \log_q\left(\frac{p}{q}\right)$$
because we have: $\log_q(S_0) = \log_q(S_m)$   (as $S_0 = S_m$)

And $\frac{\log_q\left(1+\frac{\Sigma_m}{S_0}\right)}{m} > 0$ (because $m \geq 1$ and $\Sigma_m > 0$ as the perturbation is positive)

Then, from (Lemma 3.1) we have :
$$\lim_{n\to\infty} \frac{e_n}{n} > \log_q\left(\frac{p}{q}\right) \quad \text{and} \quad \lim_{n\to\infty} \Sigma_n = +\infty$$

b) If $(S_n)_{n \geq 0}$ is unbounded then $(\lfloor S_n \rfloor)_{n \geq 0}$ is also unbounded and we have necessarily:

$$\lim_{n \to \infty} \lfloor S_n \rfloor = +\infty$$

( because again of "Pigeonhole principle", indeed if not $(\lfloor S_n \rfloor)_{n \geq 0}$ will be periodical)

We have also of course:

$$\lim_{n \to \infty} S_n = +\infty$$

It remain to show that: $\limsup_{n \to \infty} \frac{e_n}{n} \leq \log_q \left(\frac{p}{q}\right)$. Again we pass $S_n$ to logarithm:

$$\frac{e_n}{n} = \log_q \left(\frac{p}{q}\right) + \frac{\log_q\left(1+\frac{\Sigma_n}{S_0}\right)}{n} + \frac{\log_q(S_0) - 1}{n} \frac{q(S_n)}{n}$$

As $\lim_{n \to \infty} S_n = +\infty$ we have, for n large enough, $\frac{\log_q(S_0) - \log_q(S_n)}{n} \leq 0$. So it is sufficient to just prove that:

$$\lim_{n \to +\infty} \frac{\log_q\left(1+\frac{\Sigma_n}{S_0}\right)}{n} = 0$$

For that we will show that :

$$\log_q\left(1+\frac{\Sigma_n}{S_0}\right) = \sum_{j=0}^{n-1} \log_q\left(1+\frac{r_n}{S_n}\right)$$

Indeed for n ≥ 0 we have:

$$\log_q(S_0 + \Sigma_{n+1}) - \log_q(S_0 + \Sigma_n) = \log_q\left(\frac{S_0+\Sigma_{n+1}}{S_0+\Sigma_n}\right)$$

$$= \log_q\left(\frac{S_0+\Sigma_n+\frac{q^{n+e_n}}{p^n}r_n}{S_0+\Sigma_n}\right)$$

$$= \log_q\left(1 + \frac{\frac{q^{n+e_n}}{p^n}r_n}{S_0+\Sigma_n}\right)$$

$$= \log_q\left(1 + \frac{r_n}{S_n}\right)$$

So, as :

$$\log_q\left(1+\frac{\Sigma_n}{S_0}\right) = \log_q(S_0 + \Sigma_n) - \log_q(S_0)$$

$$= \sum_{j=0}^{n-1}\left(\log_q(S_0+\Sigma_{j+1}) - \log_q(S_0+\Sigma_j)\right)$$

We have:

$$\log_q\left(1+\frac{\Sigma_n}{S_0}\right) = \sum_{j=0}^{n-1} \log_q\left(1+\frac{r_n}{S_n}\right)$$

And (as $0 < r_n < 1$) :

$$\lim_{n \to \infty} S_n = +\infty \Rightarrow \lim_{n \to \infty} \log_q\left(1+\frac{r_n}{S_n}\right) = 0$$

Now by using the known Cesàro Lemma ( the convergence of the sequence of arithmetic means to the same limit) we have

$$\lim_{n \to \infty} \frac{\sum_{j=0}^{n-1} \log_q\left(1+\frac{r_j}{S_j}\right)}{n} = \lim_{n \to \infty} \log_q\left(1+\frac{r_n}{S_n}\right) = 0$$

And finally we get:

$$\lim_{n \to +\infty} \frac{\log_q\left(1+\frac{\Sigma_n}{S_0}\right)}{n} = 0 \qquad \blacksquare$$

At this stage, we are able to determine the convergence nature for the two asymptotic cases of Branch sequences.

**Lemma 3.3** *If $(S_n)_{n\geq 0}$ is a bounded, deterministic, perturbation controlled and positive Branch sequence then* $\min_{n\geq 0} S_n \leq q^2$

**Proof.** If $(S_n)_{n\geq 0}$ is a bounded then from (Lemma 3.1) and (Lemma 3.2) $(S_n)_{n\geq 0}$ is periodic (we suppose it fully periodic) of period $m \geq 1$.

Now as $S_0 = S_m$ we have:

$$\frac{p^m}{q^{m+e_m}} = \frac{S_m}{S_0 + \Sigma_m} < 1$$

and for $n = km + l$ with $k \geq 0$ and $0 \leq l < m$ we have:

$$\frac{q^{n+e_n}}{p^n} = \left(\frac{q^{m+e_m}}{p^m}\right)^k \left(\frac{q^{l+e_l}}{p^l}\right)$$

And

$$\lim_{n\to\infty} \frac{q^{n+e_n}}{p^n} = \lim_{k\to\infty} \left(\frac{q^{m+e_m}}{p^m}\right)^k \left(\frac{q^{l+e_l}}{p^l}\right) = +\infty$$

And as $(r_n)_{n\geq 0} > 0$ take a finite number of positive values because it is periodic we have also:

$$\lim_{n\to\infty} \frac{q^{n+e_n}}{p^n} r_n = +\infty$$

by (Lemma 2.3) and (Definition 2.2) of $(\omega_n)_{n\geq 0}$ let take:

$$\left(\frac{q^{h+e_h}}{p^h}\right) r_h = \max_{0\leq l < m} \left(\frac{q^{l+e_l}}{p^l} r_l\right) = \omega_h$$

Then for $n = km + h$ with $k \geq 0$ we have:

$$\left(\frac{q^{m+e_m}}{p^m}\right)^k \left(\frac{q^{h+e_h}}{p^h}\right) r_h = \max_{0\leq l < m} \left(\frac{q^{m+e_m}}{p^m}\right)^k \left(\frac{q^{l+e_l}}{p^l} r_l\right) = \omega_n$$

So by (Theorem 2.1 and (2.8.e)) we have

$$\left[\frac{\Delta_n}{q^2}\right] = 0 \text{ (which imply } \Delta_n < q^2)$$

By (Definition2.1) we have also

$$0 < S_n = C_n + \Delta_n < C_n + q^2$$

But as: $C_n = \frac{p^n}{q^{n+e_n}} \xi$ and $S_n = S_{km+h} = S_h$

We have : $S_h \leq \lim_{k \to \infty} \frac{p^n}{q^{n+e_n}} \xi + q^2 = q^2$

Which establishes that: $\min_{n \geq 0} S_n \leq q^2$

∎

**Lemma 3.4** *If $(S_n)_{n \geq 0}$ is an unbounded, deterministic, perturbation controlled and positive Branch sequence then* $\boldsymbol{limsup_{n \to \infty} \frac{e_n}{n} = log_q\left(\frac{p}{q}\right)}$

**Proof.** From (Lemma 3.2), we have seen that if $(S_n)_{n \geq 0}$ is unbounded, then:

$$\limsup_{n \to \infty} \frac{e_n}{n} \leq \log_q\left(\frac{p}{q}\right)$$

If we suppose that $\limsup_{n \to \infty} \frac{e_n}{n} < \log_q\left(\frac{p}{q}\right)$ then as $\limsup_{n \to \infty} \sqrt[n]{r_n} = 1$ we have:

$$\begin{cases} \limsup_{n \to \infty} \sqrt[n]{\frac{q^{n+e_n}}{p^n}} < 1 \\ \limsup_{n \to \infty} \sqrt[n]{\frac{q^{n+e_n}}{p^n} r_n} < 1 \end{cases}$$

So the two next series are convergent (Lemma 3.1):

$$\begin{cases} \sum_{n \geq 0} \frac{q^{n+e_n}}{p^n} r_n = \lim_{n \to \infty} \Sigma_n = \Sigma < +\infty \\ \sum_{n \geq 0} \frac{q^{n+e_n}}{p^n} < +\infty \end{cases}$$

So we have :

$$\begin{cases} \lim_{n \to \infty} \frac{q^{n+e_n}}{p^n} r_n = 0 \\ \lim_{n \to \infty} \frac{q^{n+e_n}}{p^n} = 0 \end{cases}$$

Now as:

$$\begin{cases} \lim_{n \to \infty} \frac{q^{n+e_n}}{p^n} r_n = 0 \\ \omega_n = \max_{0 \leq j \leq n-1}\left(r_j \frac{q^{j+e_j}}{p^j}\right) \end{cases}$$

the sequence $(\omega_n)_{n \geq 0}$ will be invariant for n large enough i.e.:

$$\exists K \in \mathbb{N}, \forall n \geq K, \omega_n = \omega < +\infty$$

So from (Corollary 2.2) and for n large enough we have:

$$\omega < \Sigma_n < \omega + \frac{q^{n+e_n}}{p^n} q^2$$

and hence $(as \, (\Sigma_n)_{n \geq 0}$ is a growing sequence)

$$\omega < \Sigma \leq \omega + \lim_{n \to \infty} \frac{q^{n+e_n}}{p^n} q^2 = \omega$$

Which is obviously absurd. So in conclusion we have necessarily

$$\limsup_{n \to \infty} \frac{e_n}{n} = \log_q \left(\frac{p}{q}\right) \qquad \blacksquare$$

Now we will state the better result we have obtained for the convergence of branch sequences. By (Remark 1.2) we have $min_{n \geq 0} S_n \geq 1$, so we will, *without loss of generality,* choose $S_0 = min_{n \geq 0} S_n \geq 1$.

**Theorem 3.1 (Convergence of Branch sequences)**
Let's $(S_n)_{n \geq 0}$ be a deterministic and positive perturbation Branch sequence. We suppose that:

- $S_0 = min_{n \geq 0} S_n \geq 1$ (3.1.a)
- $r_n = \frac{c}{q^{g_n}}, \forall n \geq 0$ with $c > 0$ a positive constant (3.1.b)
  (perturbation of Syracuse type)
- $\liminf_{n \to \infty} (g_n) = 0$ (3.1.c)

Then $(S_n)_{n \geq 0}$ is bounded and $S_0 = min_{n \geq 0} S_n \leq q^2$

**Proof.** Firstly with the condition (3.1.b) we have:

$$\limsup_{n \to \infty} \sqrt[n]{r_n} = \limsup_{n \to \infty} \sqrt[n]{\frac{c}{q^{g_n}}}$$

$$= \frac{1}{\liminf_{n \to \infty} \left(q^{\frac{g_n}{n}}\right)}$$

$$= \frac{1}{\left(q^{\liminf_{n \to \infty} \frac{g_n}{n}}\right)}$$

As $g_n \geq 0$ is a positive integer then the possible values for $\frac{g_n}{n}$ are 0 or integers that are multiples of $\frac{1}{n}$. The relation $\liminf_{n \to \infty}(g_n) = 0$ imply that there are infinite number of values $\frac{g_n}{n} = 0$ and then :

$$\liminf_{n \to \infty} \frac{g_n}{n} = 0 \Leftrightarrow \liminf_{n \to \infty} \left(q^{\frac{g_n}{n}}\right) = 1$$

So we have then a controlled perturbation with:

$$\limsup_{n \to \infty} \sqrt[n]{r_n} = 1 \ .$$

If $(S_n)_{n \geq 0}$ is bounded we have seen in (Lemma 3.3) that $min_{n \geq 0} S_n \leq q^2$.

So we suppose $(S_n)_{n \geq 0}$ unbounded then from (Lemma 3.2) and (Lemma 3.4) we have :

- $\lim_{n \to \infty} S_n = +\infty$

- $limsup_{n\to\infty} \left(\frac{e_n}{n}\right) = log_q\left(\frac{p}{q}\right)$

At this stage we have three possibilities for $\lim_{n\to\infty} \omega_n$ and the sequence $\left(\frac{q^{n+e_n}}{p^n}\right)_{n\geq 0}$ :

a) Case 1: $\lim_{n\to\infty} \omega_n = +\infty$.

   Then we have $limsup_{n\to\infty}\left(\frac{q^{n+e_n}}{p^n}\right) = +\infty$ because:

   $\lim_{n\to\infty} \omega_n = +\infty \Rightarrow limsup_{n\to\infty}\left(\frac{q^{n+e_n}}{p^n} r_n\right) = +\infty$

   And as $r_n = \frac{c}{q^{g_n}}$ we have

$$limsup_{n\to\infty}\left(\frac{q^{n+e_n}}{p^n} r_n\right) = limsup_{n\to\infty}\left(\frac{q^{n-1+e_{n-1}}}{p^{n-1}} \frac{qc}{p}\right)$$
$$\Leftrightarrow limsup_{n\to\infty}\left(\frac{q^{n+e_n}}{p^n}\right) = +\infty \quad (3.2.a)$$

So $\lim_{n\to\infty} \omega_n = +\infty$ imply that there exists a subsequence $(S_{n_k})_{k\geq 0}$ of $(S_n)_{n\geq 0}$ such that:

$$\begin{cases} \forall k \geq 0 \; \omega_{n_k} = \frac{q^{n_k+e_{n_k}}}{p^{n_k}} r_{n_k} \\ \lim_{n\to\infty}\left(\frac{q^{n_k+e_{n_k}}}{p^{n_k}} r_{n_k}\right) = +\infty \end{cases} \quad (3.2.b)$$

An then:

$$\begin{cases} \left\lfloor\frac{\Delta_{n_k}}{q^2}\right\rfloor = 0 \;\left(\Rightarrow \Delta_{n_k} < q^2\right) \quad \text{(Lemma 2.3 (2.6))} \\ \lim_{k\to\infty} \frac{p^{n_k}}{q^{n_k+e_{n_k}}} = 0 \quad \text{(because of (3.2. a) (3.2. b))} \end{cases}$$

So by (Corollary2.2) we have $\forall k \geq 0$:

$$0 < S_{n_k} = C_{n_k} + \Delta_{n_k} < C_{n_k} + q^2$$

And as $S_0 = \min_{n\geq 0} S_n$:

$$S_0 \leq S_{n_k} < \frac{p^{n_k}}{q^{n_k+e_{n_k}}} S_0 + q^2 \xrightarrow[k\to\infty]{} q^2$$

$$\Rightarrow S_0 \leq q^2$$

b) Case 2: $\lim_{n\to\infty} \omega_n = \sup_{n\geq 0} \omega_n = \omega < +\infty$ and $liminf_{n\to\infty}\left(\frac{q^{n+e_n}}{p^n}\right) = 0$.

   Then we have:

   $\Sigma_n - \omega_n = \sum_{j=0}^{n-1} r_j \frac{q^{j+e_j}}{p^j} - max_{0\leq j\leq n-1}\left(r_j \frac{q^{j+e_j}}{p^j}\right)$

   And we have a subsequence $(S_{n_l})_{l\geq 0}$ of $(S_n)_{n\geq 0}$ with :

   $\forall l \geq 0, \omega_{n_l} = max_{0\leq j\leq n_l-1}\left(r_j \frac{q^{j+e_j}}{p^j}\right) = r_{n_l}\frac{q^{n_l+e_{n_l}}}{p^{n_l}}$

   So $\Sigma_{n_l} - \omega_{n_l} = \Sigma_{n_l-1}$ which is a growing sequence and by passing to the limit we have:

$$\lim_{l\to\infty}\Sigma_{n_l} - \omega = \lim_{l\to\infty}\Sigma_{n_l}$$

As $\omega > 0$ (because the perturbation is positive (Definition 1.4)) we have then perforce:

$\lim_{n\to\infty}\Sigma_n = +\infty$ and for n sufficiently large we have $\Sigma_n > \omega$

Now we have $liminf_{n\to\infty}\left(\frac{q^{n+e_n}}{p^n}\right) = 0$ so again there exists a subsequence $(S_{n_k})_{k\geq 0}$ of $(S_n)_{n\geq 0}$ such that: $lim_{k\to\infty}\left(\frac{q^{n_k+e_{n_k}}}{p^{n_k}}\right) = 0$

From (Corollary 2.2) we have $\forall n \geq 0, \; \omega_n < \Sigma_n < \omega_n + \frac{q^{n+e_n}}{p^n}q^2$

So we have finally for k sufficiently large:

$$\omega < \Sigma_{n_k} < \omega + \frac{q^{n_k+e_{n_k}}}{p^{n_k}}q^2$$

And $(as\;(\Sigma_{n_k})_{k\geq 0}$ is a growing sequence$)$

$$\omega < \lim_{k\to\infty}\Sigma_{n_k} \leq \omega + \lim_{n\to\infty}\frac{q^{n+e_n}}{p^n}q^2 = \omega$$

This is obviously absurd.

c) Case 3: $\lim_{n\to\infty}\omega_n = \sup_{n\geq 0}\omega_n = \omega < +\infty$ and $liminf_{n\to\infty}\left(\frac{q^{n+e_n}}{p^n}\right) > 0$.

By (Definition 2.2) of $(\omega_n)_{n\geq 0}$ we have

$$limsup_{n\to\infty}\left(\frac{q^{n+e_n}}{p^n}r_n\right) < +\infty$$

And as $r_n = \frac{c}{q^{g_n}}$ we have

$$limsup_{n\to\infty}\left(\frac{q^{n+e_n}}{p^n}r_n\right) = limsup_{n\to\infty}\left(\frac{q^{n-1+e_{n-1}}}{p^{n-1}}\frac{qc}{p}\right)$$

So we have also

$$limsup_{n\to\infty}\left(\frac{q^{n+e_n}}{p^n}\right) < +\infty$$

So $\left(\frac{q^{n+e_n}}{p^n}\right)_{n\geq 0}$ is a bounded sequence in a strictly subinterval of $]0,+\infty[$ (and this is also the case for $\left(\frac{p^n}{q^{n+e_n}}\right)_{n\geq 0}$). Now since we have from (Corollary 2.2):

$$\forall n \geq 0, \omega_n < \Sigma_n < \omega_n + \frac{q^{n+e_n}}{p^n}q^2$$

The sequence $(\Sigma_n)_{n\geq 0}$ is also bounded and consequently as:

$$S_n = \left(\frac{p^n}{q^{n+e_n}}\right)(\xi + \Sigma_n)$$

The sequence $(S_n)_{n\geq 0}$ will be also bounded and by (Lemma 3.3) we have:

$S_0 = min_{n\geq 0} S_n \leq q^2$

■

## 4. Application to Syracuse sequences

The Syracuse problem (also called the Collatz problem or $(3x+1)$ problem) may be stated in a variety of ways [1], for example by defining the next sequence $(T_n)_{n\geq 0}$.

**Definition 4.1** A Syracuse sequence $(T_n)_{n\geq 0}$ is a sequence of positive integers given by:

- $T_0 > 0$
- for $n \geq 0$ we have:
$$\begin{cases} T_{n+1} = \frac{3T_n+1}{2} & \text{if } T_n \equiv 1 \bmod(2) \\ T_{n+1} = \frac{T_n}{2} & \text{if } T_n \equiv 0 \bmod(2) \end{cases} \tag{4.1}$$

□

We have then the next conjecture

**Conjecture 4.1 (Syracuse Sequence Convergence)** *For any positive integer number $T_0 > 0$, the sequence $(T_n)_{n\geq 0}$ ends with the trivial cycle $(1,2)$.*

□

A formal proof is lacking so far in spite of various approaches to the problem [2] [3] [4] [5] [6] [7] [8] [9] [20].

Another equivalent formulation where the sequence $(T_n)_{n\geq 0}$ is reduced to only odd numbers is defined as follows.

**Definition 4.2** We consider the sequence $(W_n)_{n\geq 0}$ of odd numbers, called also Syracuse sequence and defined by:

- $W_0$ is a given odd positive integer and $h_0 = 0$
- for $n \geq 0$, $W_{n+1} = \frac{3W_n+1}{2^{1+h_{n+1}}}$ with $h_{n+1} = v_2\left(\frac{3W_n+1}{2}\right)$ (4.2)

($v_2$ is the 2-adique valuation: the exponent of the highest power of the prime number 2 that divides $\frac{3W_n+1}{2}$)

□

This definition provides, exactly as we have seen in (Lemma 1.1), a direct expression for $(W_n)_{n\geq 0}$

**Lemma 4.1** *We can have a direct expression of $W_n$ in terms of $W_0$ and $(h_k)_{k<n}$:*
$$\forall n \geq 0, \ W_n = \left(\frac{3^n}{2^{n+e'_n}}\right)(W_0 + \Sigma'_n) \tag{4.3}$$
With the sequences $(\Sigma'_n)_{n\geq 0}$ and $(e'_n)_{n\geq 0}$ are given by:
$$\begin{cases} \Sigma'_0 = 0 \\ e'_0 = h_0 = 0 \end{cases} \text{ and } \begin{cases} \Sigma'_n = \frac{1}{3}\sum_{j=0}^{n-1}\frac{2^{j+e'_j}}{3^j} \\ e'_n = \sum_{j=0}^{n} h_j \end{cases}$$

**Proof.** The same proof as in (Lemma 1.1) can be used with p=3 and q=2.

∎

Now we will introduce a Branch sequence $(S_n)_{n\geq 0}$ that can be associated to the Syracuse sequence $(W_n)_{n\geq 0}$ and with same convergence properties.

**Lemma 4.2** *We set $p = 3, q = 2$ and $g_0 = 0$, $\forall n \geq 1$ $g_n = h_{n+1}$. We define the sequence $(S_n)_{n \geq 0}$ by:*

$$\forall n \geq 0, \quad S_n = \frac{qW_n}{q^{h_{n+1}}} \qquad (4.3.a)$$

*Then $(S_n)_{n \geq 0}$ is a Branch sequence with a controlled positive perturbation sequence $(r_n)_{n \geq 0}$ given by:*

$$\forall n \geq 0, \quad r_n = \frac{q}{pq^{g_n}} \qquad (4.3.b)$$

*i.e. the sequence $(S_n)_{n \geq 0}$ verify (Definition(1.3)) and:*

$$\forall n \geq 0, \quad S_n = \left(\frac{p^n}{q^{n+e_n}}\right)(S_0 + \Sigma_n) \qquad (4.3.c)$$

*With: $e_n = \sum_{j=0}^{n} g_j$, $\Sigma_0 = 0$ and $\Sigma_n = \sum_{j=0}^{n-1} \frac{q^{j+e_j}}{p^j} r_j$*

**Proof.** Firstly, the values p=3 and q=2 verifies the (Definition 1.1) conditions about $(p, q)$:
  $p > q > 1$ *are coprime integers,*
  $p = \alpha q + \beta$, *with* $\alpha = \beta = 1$
  $\alpha + \beta = 2 = q$

Then by (Definition 4.1) we have:
$$\forall n \geq 0, \quad W_{n+1} = \frac{pW_n + 1}{q^{1+h_{n+1}}}$$

And :
$$\begin{aligned}
S_{n+1} &= \frac{qW_{n+1}}{q^{h_{n+2}}} \\
&= \frac{q}{q^{h_{n+2}}} \frac{pW_n + 1}{q^{1+h_{n+1}}} \\
&= \frac{p\left(\frac{qW_n}{q^{h_{n+1}}} + \frac{q}{pq^{h_{n+1}}}\right)}{q^{1+h_{n+2}}} \\
&= \frac{p(S_n + r_n)}{q^{1+g_{n+1}}}
\end{aligned}$$

So, as in (Lemma1.1), we have the (4.3.c) expression:
$$\forall n \geq 0, \quad S_n = \left(\frac{p^n}{q^{n+e_n}}\right)(S_0 + \Sigma_n)$$

We have also (4.3.b): the sequence $(r_n)_{n \geq 0}$ is given by :
$$\forall n \geq 0, \quad r_n = \frac{q}{pq^{g_n}}$$

Now as $W_n$ is an odd number, we have $2W_n \equiv 2 \mod(2^2)$ and we have the binary expansion of $2W_n$ begin always with a sequence of alternate zero and one 10 and finish necessarily with either a double zero 00 or a double one 11:

$$\begin{cases} 2W_n = (\cdots 001010 \cdots\cdots 010)_2 \\ or \\ 2W_n = (\cdots 11010 \cdots\cdots 010)_2 \end{cases}$$

This means that, for some positive integers $a \geq 0$ and $b > 0$ we have one of the two next cases for $2W_n$:

$$\begin{cases} 2W_n = 2\frac{4^a - 1}{3} + 4^a 2b \text{ with } b \equiv 1 \mod 8 & (4.4.a) \\ or \\ 2W_n = 2\frac{4^a - 1}{3} + 4^a b \text{ with } b \equiv 6 \mod 8 & (4.4.b) \end{cases}$$

Or as shown in the next table:

$$\left(\underbrace{\cdots 001\ 0\ \overbrace{10\cdots\cdots 010}^{2^{\frac{4^a-1}{3}}}}_{2W_n}\right)_2 \quad or \quad \left(\underbrace{\cdots \overbrace{110}^{b}\ \overbrace{101\cdots\cdots 010}^{2^{\frac{4^a-1}{3}}}}_{2W_n}\right)_2 \qquad (4.4.c)$$

We have hence

$$\begin{cases} \frac{3W_n+1}{4^a} = 1 + 3b \equiv 4 \bmod 8 \\ or \\ \frac{3W_n+1}{4^a} = 1 + 3\frac{b}{2} \equiv 2 \bmod 8 \end{cases}$$

So as $g_n = h_{n+1} = \vartheta_q\left(\frac{3W_n+1}{2}\right)$ and $S_{n+1} = \frac{2W_n}{2^{h_{n+1}}}$ we have:

$$\begin{cases} h_{n+1} = 2a + 1 \text{ and } \lfloor S_n \rfloor = 2b \text{ for the case } (4.4.a) \\ or \\ h_{n+1} = 2a \text{ and } \lfloor S_n \rfloor = b \text{ for the case } (4.4.b) \end{cases}$$

In the two cases, conditions (1.3.c) and (1.3.d) are then fulfilled for the sequence $(S_n)_{n\geq 0}$ as shown in the next table:

$$\left(\underbrace{\overbrace{\cdots 0010}^{\lfloor S_n \rfloor}\ \overbrace{10\cdots\cdots\cdot 01010}^{2^{h_{n+1}}\{S_n\}=2^{\frac{2^{h_{n+1}-1}-1}{3}}}}_{2W_n}\right)_2 \quad or \quad \left(\underbrace{\overbrace{\cdots 110}^{\lfloor S_n \rfloor}\ \overbrace{101\cdots\cdots\cdot 01010}^{2^{h_{n+1}}\{S_n\}=2^{\frac{2^{h_{n+1}-1}-1}{3}}}}_{2W_n}\right)_2 \qquad (4.4.d)$$

From (4.4.d) we have also:

$$\begin{cases} \{S_n\} = \frac{1}{3} - \frac{2}{3 \cdot 2^{h_{n+1}}} & \text{for the case } (4.4.a) \text{ (in this case } h_{n+1} \geq 1) \\ or \\ \{S_n\} = \frac{2}{3} - \frac{2}{3 \cdot 2^{h_{n+1}}} & \text{for the case } (4.4.b) \end{cases} \qquad (4.4.e)$$

So as $r_n = \frac{2}{3 \cdot 2^{h_{n+1}}}$ we have in the two cases:
$$0 \leq \{S_n\} + r_n < 1 \quad \text{i.e.} \quad \lfloor S_n + r_n \rfloor = \lfloor S_n \rfloor$$

∎

**Theorem 4.1** *(Convergence of Syracuse sequences) A Syracuse sequences $(W_n)_{n\geq 0}$ can be associated to deterministic Branch sequences $(S_n)_{n\geq 0}$ with perturbation of Syracuse type :*

$$\begin{cases} p = 3, q = 2 \\ r_{n=\frac{q}{p*q^{g_n}}} \end{cases}$$

*And this confirms the Conjecture 4.1.*

**Proof.** (Lemma 4.2) establishes that a Syracuse sequence $(W_n)_{n\geq 0}$ (Definition 4.1) can be associated to a Branch sequence $(S_n)_{n\geq 0}$ (Definition 1.3) with a perturbation sequence $(r_n)_{n\geq 0}$ given by:  $\forall n \geq 0, \ r_n = \frac{q}{pq^{g_n}}$.

We have also:

$$W_{n+1} = \frac{3W_n+1}{2^{1+h_{n+1}}} \Leftrightarrow 2W_{n+1} = \frac{3S_n}{2} + \frac{1}{2^{h_{n+1}}}$$

So from (Lemma 4.2) and (4.4.e)

$$\begin{cases} 2W_{n+1} = \dfrac{3\lfloor S_n \rfloor}{2} + \dfrac{3\{S_n\}}{2} + \dfrac{1}{2^{h_{n+1}}} = \dfrac{3\lfloor S_n \rfloor + 1}{2} & \text{for the case (4.4.a)} \\ \quad or \\ 2W_{n+1} = \dfrac{3\lfloor S_n \rfloor}{2} + \dfrac{3\{S_n\}}{2} + \dfrac{1}{2^{h_{n+1}}} = \dfrac{3\lfloor S_n \rfloor + 1}{2} & \text{for the case (4.4.b)} \end{cases}$$

So for the two cases we have $\quad 2W_{n+1} = \dfrac{3\lfloor S_n \rfloor + 1}{2}$ \hfill (4.5.a)

This implies that $(S_n)_{n \geq 0}$ is deterministic because:

$$\lfloor S_n \rfloor = \lfloor S_m \rfloor \Rightarrow W_{n+1} = W_{m+1} \quad \text{(from (4.5.a))}$$

$$\Rightarrow \begin{cases} S_{n+1} = S_{m+1} \\ r_{n+1} = r_{m+1} \end{cases} \quad (from\ (4.3.a))$$

Now the convergence result of (Theorem 3.1) is conditioned by :

$$\liminf_{n \to \infty}(g_n) = 0$$

But if we have $\liminf\limits_{n \to \infty}(g_n) > 0$ then the sequence $(S_n)_{n \geq 0}$ is necessarily bounded because in this case and for n sufficiently large we have $g_n \geq 1$ then

$$\begin{aligned} S_{n+1} &= \frac{qW_{n+1}}{q^{g_{n+1}}} \\ &= \frac{q}{q^{g_{n+1}}} \frac{pW_n + 1}{q^{1+g_n}} \\ &= \frac{q}{q^{g_n}} \frac{pW_n + 1}{q^{1+g_{n+1}}} \\ &< \frac{qW_n}{q^{g_{n+1}}} \quad \text{if } W_n > 1 \\ &= S_n \end{aligned}$$

and $(S_n)_{n \geq 0}$ is a decreasing sequence and then necessarily $min_{n \geq 0} S_n \leq q^2$ .

So finally:

- Either $\exists n \geq 0, W_n = 1$ and the Syracuse sequence finishes by the trivial cycle and of course $min_{n \geq 0} S_n \leq q^2$
- Or $\forall n \geq 0, W_n > 1$ and we have the convergence result of (Theorem 3.1) and the sequence $(S_n)_{n \geq 0}$ is again bounded and $min_{n \geq 0} S_n \leq q^2$

In conclusion we have always $min_{n \geq 0} S_n \leq 4$ and as
$$2W_{n+1} = \frac{3\lfloor S_n \rfloor + 1}{2}$$
We have
$$min_{n \geq 0} W_n \leq 3$$
This obviously confirms the Syracuse conjecture

∎

## 5. Syracuse sequences cellular automata

We will finish this paper by presenting a visualization method for Branch Syracuse sequences $(S_n)_{n\geq 0}$ as cellular automaton (CA). This will help to understand how the Branch sequence proof works and in the same time justifies the choice of the name "Branch" for such sequences.

We begin by considering the binary expansions of the Syracuse sequence $(W_n)_{n\geq 0}$:

$$W_n = \sum_{j\in N} b_{(n,j+e_n)} 2^j = (\cdots b_{(n,2+e_n)} b_{(n,1+e_n)} b_{(n,0+e_n)})_2 \, , \qquad for\ n \geq 0,$$

Then binary coefficients $\left(b_{(n,j+e_n)}\right)_{n,j\geq 0}$ form, by their indexes, a subset of the lattice $Z^2$ and $e_n$ fixes the beginning position of integer part of $W_n$ in the lattice (the lattice is completed by zero values). Then $\left(b_{(n,j)}\right)_{n,j\geq 0}$ can be seen as a dynamical system, more precisely as a binary one dimensional cellular automaton (CA) evolving with the index $n$ [22]; the automaton is initialized with the binary expansion of $W_0$, and iteratively when the line $n$ is calculated we use the next transition rules (as defined in (2.3.c)) to generate the line $(n + 1)$:

$$\begin{cases} b_{(n+1,j)} = \left(b_{(n,j+1)} + b_{(n,j)} + \delta_{(n,j)}\right) mod(2) \\ \delta_{(n,j+1)} = \left\lfloor \frac{(b_{(n,j+1)}+b_{(n,j)}+\delta_{(n,j)})}{2} \right\rfloor \end{cases} \qquad (4.6.a)$$

Each cell, at position $(n,j)$ is represented by the couple $\left(b_{(n,j)}, \delta_{(n,j)}\right)$ with $\left(\delta_{(n,j)}\right)_{n\geq 0}$ the carry digits of the binary arithmetic addition: $2 * W_n + (W_n + 1) = 2^{1+g_{n+1}} W_{n+1}$

$$\overbrace{\cdots\cdots b_{(n,j+1)}{}^{\delta_{(n,j+1)}} b_{(n,j)}{}^{\delta_{(n,j)}} \cdots\cdots b_{(n,e_n)} \overset{1}{0} \cdots\cdots\cdots 0}^{\text{n line: } 2^{1+g_{n-1}} W_n}$$

$$\underbrace{\cdots\cdots\cdots\cdots\cdots\cdots\cdots b_{(n+1,j)} \cdots\cdots b_{(n+1,e_n)}\, 0 \cdots\cdots\cdots 0}_{\text{(n+1) line: } 2^{1+g_n} W_{n+1} = 2*W_n + (W_n+1)}$$

Juxtaposition of n line with (n+1) line in cellular automaton

So we have, $for\ j \geq e_n\ and\ with\ \delta_{(n,e_n)} = 1$, the next transition equation from line $n$ to line $(n + 1)$:

$$b_{(n,j+1)} + b_{(n,j)} + \delta_{(n,j)} = b_{(n+1,j)} + 2 * \delta_{(n,j+1)} \qquad (3.5.b)$$

To generate the Syracuse automaton line (n+1), we begin arithmetic calculations at position $e_n$ of line n and $\delta_{(n,e_n)}$ is initialized to 1. Note that divisions by 2 correspond to position shifts in the cellular automata lines. We can, by the same transition rules, generate a cellular automata associated to the sequence of rational powers $(W_0(3/2)^n)_{n\geq 0}$ but this time without any carry digits initializations. We see then ($Figure.1.1$) that the Syracuse automaton is perturbation (by the initialization of carry digits at each line) of the (3/2) rational powers automaton. We note that we have the same CA if we consider the sequence $(S_n)_{n\geq 0}$ instead of $(W_n)_{n\geq 0}$ because they have the same binary digits in their binary expansion.

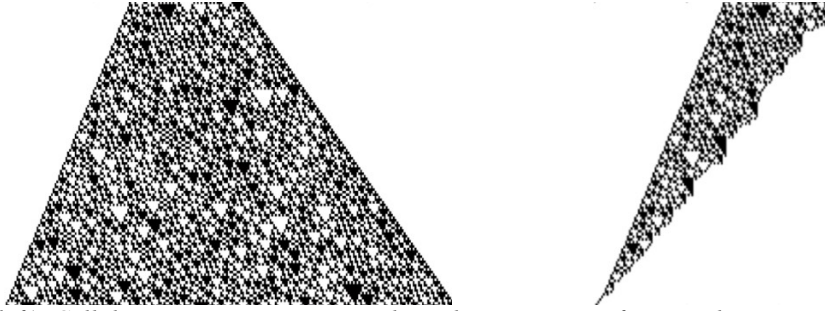

***Figure 5.1:*** *(left) Cellular automaton associated to the sequence of rational powers $(C_n = S_0(3/2)^n)_{n\geq 0}$; (right) Cellular automaton associated to Syracuse sequence $(S_n)_{n\geq 0}$. The two automata are generated with the same initialization $S_0$. The Syracuse Automaton can be seen as a perturbation of the rational powers sequence automaton by carry digit initialization at each line; we can see also the convergence of the Syracuse cellular automaton, to the cycle (1,2).*

An uncommon way of expressing binary expansions is to use the binary gray coding [18] [19]. It is mainly used in electronics, because it allows the modification of only one bit at a time when a number is increased by one. If we express each line of a (3/2) cellular automaton in the binary gray code instead of the standard binary expansion, we obtain a new cellular automaton composed of the cells $(g_{(n,j)})_{n,j\geq 0} \subset Z^2$. Switching from standard binary coding $(b_{(n,j)})_{n,j\geq 0}$ to gray coding is done very easily by:

$$g_{(n,j)} = b_{(n,j+1)} \oplus b_{(n,j)} \qquad (3.5.c)$$

where $\oplus$ is the exclusive disjunction (XOR operator) *(Figure 5.2)*.

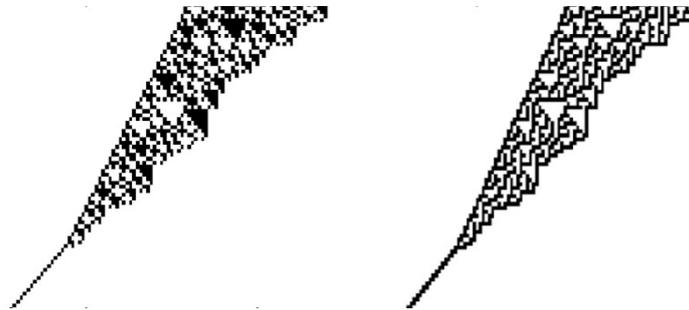

***Figure 5.2:*** *Cellular automaton associated to Syracuse sequence $(S_n)_{n\geq 0}$ expressed (left) in binary expansion and (right) in gray code. The gray code automaton brings out a characteristic tree structure with branches bifurcating only upwards. The branches are coalescing towards the automaton left edge to form a single trunk.*

When we cross the Gray code cellular automata from top to bottom the cells with values 1 extend only either to the left or to the bottom. This clearly shows that the automata have a tree structure, with branches bifurcating only upwards, and justifies the name of "Branch sequences". We have, in a way, proved with Branch sequences that as long as two branches are separated, perturbations (more precisely controlled perturbations) that happen leftward the left branch do not affect the right one. So perturbations progress to the left only through coalescing branches and this explains the perturbation confinement in the in the CA right border (Figure 5.3).

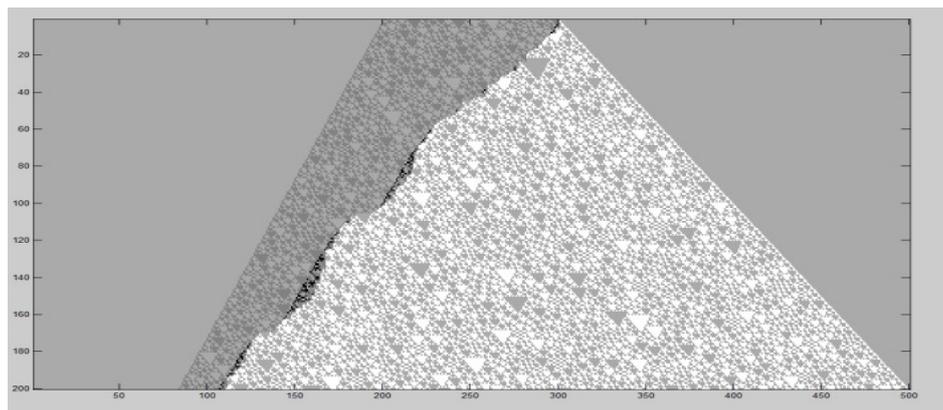

***Figure 5.3:*** *Superposition of the Cellular Automata associated to rational powers sequence $(C_n)_{n\geq 0}$ and its associated Syracuse sequence $(S_n = C_n + \Delta_n)_{n\geq 0}$: (Left) we have the identical part between $(S_n)_{n\geq 0}$ and $(C_n)_{n\geq 0}$. (Right) the fractional part of the perturbation $(\{\Delta_n\})_{n\geq 0}$. (Between) the perturbation integer part $(\lfloor \Delta_n \rfloor)_{n\geq 0}$. We see then the clearly the perturbation confinement phenomena into the Syracuse CA right border.*

To sum up, all the branches in the gray code CA seem (Figure 5.2) to converge towards the left edge of the automaton to form a single trunk. This is, somehow, another formulation of the Syracuse conjecture:

"All Gray code automata of Syracuse sequences have a "tree" structure with coalescing "branches.""

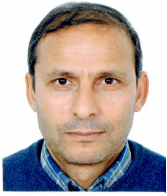

Hassan Douzi is a research professor in mathematics at Ibn Zohr University, Agadir, Morocco. He received his PhD in "Mathematics for decision" from the University of Paris IX Dauphine 1992. He published several research papers, essentially in Applied mathematics related to image processing, but also incidentally preprint papers in Number theory and Combinatorics.